\documentclass[lfeqn,12pt]{article}
\usepackage{amssymb,bbm,pifont,mathtools,mathabx,mathrsfs,
}
\usepackage{minitoc}
 
\usepackage[nottoc]{tocbibind}

\usepackage{comment}

\usepackage{multirow}
\usepackage[multiple]{footmisc}
\usepackage{xcolor}
\usepackage[colorlinks=true]{hyperref}
\usepackage{marginnote}

\newcommand{\Giu}{{\bigskip\noindent}}

\newcommand{\noi}{{\noindent}}

\newtheorem{theorem}{Theorem}
\newtheorem{definition}[theorem]{Definition}
\newtheorem{proposition}[theorem]{Proposition}
\newtheorem{lemma}[theorem]{Lemma}
\newtheorem{remark}[theorem]{Remark}
\newtheorem{conjecture}[theorem]{Conjecture}
\newtheorem{question}{Question}
\newtheorem{sublemma}[theorem]{Sublemma}
\newtheorem{corollary}[theorem]{Corollary}
\newtheorem{assumption}[theorem]{Assumption}
\newtheorem{notationalrem}[theorem]{Notational Remark}

\newtheorem{tools}[subsection]{$\negsp\negsp$}
\newcommand\asm[1]{ \begin{assumption}\label{#1} }
\newcommand\easm{ \end{assumption} }
\newcommand\dfn[1]{ \begin{definition}\label{#1} }
\newcommand\dfntwo[2]{ \begin{definition}[#2]\label{#1} }
\newcommand\edfn{ \end{definition} }
\newcommand\rem[1]{ \begin{remark}\label{#1} \small \rm}
\newcommand\remtwo[2]{ \begin{remark}[#2]\label{#1} \rm}
\newcommand\erem{ \end{remark} }
\newcommand\thm[1]{ \begin{theorem}\label{#1}}
\newcommand\thmtwo[2]{ \begin{theorem}[#2]\label{#1}}
\newcommand\ethm{ \end{theorem} }
\newcommand\pro[1]{ \begin{proposition}\label{#1}}       
\newcommand\protwo[2]{ \begin{proposition}[#2]\label{#1}}
\newcommand\epro{ \end{proposition} }
\newcommand\conj[1]{ \begin{conjecture}\label{#1}}       
\newcommand\conjtwo[2]{ \begin{conjecture}[#2]\label{#1}}
\newcommand\econj{ \end{conjecture} }
\newcommand\quest[1]{ \begin{question}\label{#1}}       
\newcommand\questtwo[2]{ \begin{question}[#2]\label{#1}}
\newcommand\equest{ \end{question} }
\newcommand\lem[1]{ \begin{lemma}\label{#1}}
\newcommand\lemtwo[2]{ \begin{lemma}[#2]\label{#1}}
\newcommand\elem{ \end{lemma} }
\newcommand\sublem[1]{ \begin{sublemma}\label{#1}}
\newcommand\sublemtwo[2]{ \begin{sublemma}[#2]\label{#1}}
\newcommand\esublem{ \end{sublemma} }
\newcommand\cor[1]{ \begin{corollary}\label{#1}}
\newcommand\cortwo[2]{ \begin{corollary}[#2]\label{#1}}
\newcommand\ecor{ \end{corollary} }
\newcommand\notrem[1]{ \begin{notationalrem}\label{#1} \sl}
\newcommand\enotrem{ \end{notationalrem} }

\newcommand\equ[1]{{\rm (\ref{#1})}}
\newcommand\beq[1]{ \begin{equation}\label{#1} }
\newcommand{\eeq}{ \end{equation} }

\newcommand\beqa[1]{ \begin{eqnarray} \label{#1}}
\newcommand{\eeqa}{ \end{eqnarray} }
\newcommand{\beqano}{ \begin{eqnarray*} }
\newcommand{\eeqano}{ \end{eqnarray*} }

\newcommand{\proof}{\par\medskip\noindent{\bf Proof\ }}

\newcommand\warning[1]{\textcolor{red}{{#1}}}

\newcommand{\ie}{{\it i.e.\  }}

\newcommand{\eg}{{\it e.g. \ }}

\newcommand{\wrt}{{\it w.r.t\ }}

\newcommand{\sss}{{\it iff\ }}

\newcommand{\qed}{\hskip.5truecm
            \vrule width 1.7truemm height 3.5truemm depth 0.truemm
            \par\Giu}

\newcommand\su[1]{ \frac{1}{ {#1}} }

\newcommand\jac{ {\, \rm J }}
\newcommand\len{ {\, \rm length }}
\newcommand\am{ {\, \rm am }}
\newcommand\ecn{ {\, \rm cn }}
\newcommand\esn{ {\, \rm sn }}

\newcommand{\rot}{ {\rm \, rot \, }}

\newcommand{\dpr}{ {\partial}   }

\newcommand\eqby[1]{\stackrel{\equ{#1}}{=}}

\newcommand{\negsp}{\hspace{-.04truecm}}

\renewcommand{\a }{ {\alpha}   }
\renewcommand{\b}{ {\beta}   }
\newcommand{\g}{ {\gamma}   }
\newcommand{\G}{ {\Gamma}   }
\renewcommand{\d}{ {\delta}   }

\newcommand{\vae }{ {\varepsilon}   }

\renewcommand{\th }{ {\theta}   }

\renewcommand{\l}{ {\lambda}   }

\newcommand{\x }{ {\xi}   }

\newcommand{\s}{ {\sigma}   }

\newcommand{\f}{ {\varphi}   }

\renewcommand{\o}{ {\omega}   }
\renewcommand{\O}{ {\Omega}   }
\newcommand{\torus}{ {\mathbb{ T}}   }
\renewcommand{\natural}{ {\mathbb{ N}}   }
\newcommand{\real}{ {\mathbb{ R}}   }
\newcommand{\integer}{ {\mathbb{ Z}}   }

\newcommand{\rational}{ {\mathbb{ Q}}  }

\newcommand{\rn}{ {\real^2}   }

\newcommand{\zn}{ {\integer^d }   }

\font\teneufm=eufm10
\font\seveneufm=eufm7
\font\fiveeufm=eufm5
\newfam\eufmfam
\textfont\eufmfam=\teneufm
\scriptfont\eufmfam=\seveneufm
\scriptscriptfont\eufmfam=\fiveeufm

\newcommand\appA[1]{\section{#1}\label{app:A}
\renewcommand{\theequation}{A.\arabic{equation}}
           \setcounter{equation}{0}
\renewcommand{\thetheorem}{A.\arabic{theorem}}
           \setcounter{theorem}{0}
                  }
\newcommand\appB[1]{\section{#1}\label{app:B}
\renewcommand{\theequation}{B.\arabic{equation}}
           \setcounter{equation}{0}
\renewcommand{\thetheorem}{B.\arabic{theorem}}
           \setcounter{theorem}{0}           
           }
\newcommand\appC[1]{\section{#1}\label{app:C}
\renewcommand{\theequation}{D.\arabic{equation}}
           \setcounter{equation}{0}
\renewcommand{\thetheorem}{D.\arabic{theorem}}
           \setcounter{theorem}{0}           
           }

\newcommand{\id}{{\mathtt {id}}}

\renewcommand\subset{\subseteq}

\begin{document}


\date{\small \today}

\title{
{\bf  
On some invariants of Birkhoff billiards under conjugacy
}\\
}
\author{
V. Kaloshin\footnote{Vadim Kaloshin: University of Maryland, College Park, MD, USA \& 
Institute of Science and Technology Austria, Am Campus~1, 3400 Klosterneuburg, Austria. 
\emph{vadim.kaloshin@gmail.com}}\ , 
\
C. E. Koudjinan\footnote{Comlan Edmond Koudjinan: Institute of Science and Technology Austria (IST Austria), Am Campus~1, 3400 Klosterneuburg,
Austria. \emph{edmond.koudjinan@ist.ac.at}
} 
}
\maketitle
\begin{abstract}
In  the class of strictly convex smooth boundaries, each of which not having strip around its boundary foliated by invariant curves, we prove that the Taylor coefficients of the "normalized" Mather's $\beta$--function  
are invariants under $C^\infty$--conjugacies.  
In contrast, we prove that any two elliptic billiard maps are $C^0$--conjugated near their respective boundaries, and $C^\infty$--conjugated in the open cylinder, near the boundary and away from a plain passing through the center of the underlying ellipse. We also prove that 
if the billiard maps corresponding to two ellipses are topologically conjugated then the two ellipses are similar. 
\end{abstract}

\noi
\section{Introduction}

A billiard is a mathematical modeling of the dynamic of a confined massless particle without friction and reflecting elastically on the boundary (without friction): the particle moves along a straight line with constant speed till it hits the boundary, then reflects off with reflection angle equals to the angle of incidence and follows the reflected straight line. This seemingly simple and fascinating Dynamical system was introduced by G.D. Birkhoff \cite{birkhoff2020periodic} in 1920. Since then, it has captured much of attention from both Physicists and Mathematicians. We refer the reader to \cite{tabachnikov2005geometry} and references therein for more details on Billiards. It is presently a very active and popular 
subject with many challenging open problems.

\noi
The motivation of this work is the following question attributed to Victor Guillemin:
 \quest{qes1}
 Let $f$ and $g$ be smooth Birkhoff billiard maps corresponding to two
strictly convex domains $\O_f$ and $\O_g$. Assume that f and g are conjugate, i.e. there exists a homeomorphism h such that $f = h^{-1} \circ g \circ h$. What can we say about the two domains $\O_f$ and $\O_g$? Are they ``similar", that is, have they the same shape?
 \equest
 To our knowledge the only known answer to this question is in the case of circular billiards. In fact, a billiard map in a disc $\mathscr D$ enjoys the peculiar property of having the phase space completely foliated by homotopically non-trivial invariant curves.
It is an example of so--called globally integrable billiards. In terms of the geometry of the billiard domain $\mathscr D$, this reflects the existence of a smooth foliation by (smooth and
convex) caustics, i.e. (smooth and convex) curves with the property that if a trajectory is tangent to one of them, then it will remain tangent after each reflection. It is easy to check that in the circular billiard case this family consists of concentric circles. See for instance \cite[Chapter 2]{tabachnikov2005geometry}. 
The peculiarity of this property is justified by the following exciting result due to Misha Bialy: 

\ \\
\noi
{\bf Theorem (M. Bialy\cite{bialy1993convex})} If the phase space of the billiard ball map is foliated by continuous invariant curves which are not null--homotopic, then it is a circular billiard.

\ \\
\noi
In particular, the answer to Question~\ref{qes1} is affirmative whenever one of the domains is a circle, since any foliation is invariant under conjugacy.

\ \\
\noi
In general, billiards in ellipses are locally integrable: a neighborhood of their boundary is foliated by caustics. Then, Birkhoff asked 
whether there are other examples of (locally) integrable billiards. 
He actually conjectured:
\conjtwo{conj1}{Birkhoff}
Amongst all convex billiards, only those in ellipses (circles being a distinct special case) are integrable.
\econj
We are mainly concerned in this paper by invariant of Birkhoff billiards under conjugacies when both of the domains are non--integrable and with the existence of conjugacy (local and gloabal) for elliptic billiards tables (which corresponds to the integrable case according to Birkhoff conjecture). 
before summarizing the results of this work, let us introduce the notion of "normalized" Mather's $\beta$--function. 

\noi
The Mather's $\beta$--function (or Mather's minimal average action) is a key function in the celebrated Aubry--Mather Theory as it encodes several important properties of the dynamics. The Mather's $\beta$--function can be defined for any exact area preserving twist map, not
necessarily a billiard map. Roughly speaking, its associates to any fixed rotation number (not only
rational ones) the minimal average action of orbits with that rotation number
(whose existence, inside a suitable interval, is ensured by the twist condition), which we now define precisely.

\noi
For, let $\O$ be a bounded strictly convex domain in $\real^2$ with $C^r$ boundary $\dpr\O$, with $r\ge 3$.\footnote{Observe that if $\O$ is not convex, then the billiard map is not
continuous; in this article we will be interested only in strictly convex domains. Moreover, as pointed out by Halpern \cite{halpern1977strange}, if the boundary is
not at least $C^3$, then the flow might not be complete.} 
The phase space $\mathscr M$ of the billiard map consists of unit vectors $(x, v)$ whose foot points $x$ are on $\dpr\O$ and that have inward directions.  The billiard ball map $f : \mathscr M\to \mathscr M$ takes $(x, v)$ to $(x',v')$, where $x'$ represents the point where
the trajectory starting at $x$ with velocity $v$ hits the boundary $\dpr\O$ again, and $v'$ is the reflected velocity, according to the standard reflection law: angle of incidence is equal to the angle of reflection. 
\noi
Assume that the boundary $\dpr\O$ is parametrized by arc--length $s$ and let $\g : [0, |\dpr\O|]  \to \real^2$ denote such a parametrization,
where $|\dpr\O|$ denotes the length of $\dpr\O$. Let $\th$ be the angle between $v$ and the positive tangent to $\dpr\O$ at $x$. Hence, $\mathscr M$ can be identified with the annulus $\mathbb{A} \coloneqq [0, |\dpr\O|] \times (0, \pi)$ and the billiard map $f$ can be described as
$$
f\colon \mathbb{A}\ni (x,\th)\longmapsto (x',\th')\in \mathbb{A}.
$$
Moreover, $f$ can be extended (smoothly) to $\bar{\mathbb{A}}=[0, |\dpr\O|] \times [0, \pi]$ by setting $f(s,\th_0)=(s,\th_0)$, for all $s\in[0, |\dpr\O|]$, whenever $\th_0\in \{0,\pi\}$. Denote the Euclidean distance between two points of $\dpr\O$ by
$$
d(s,s')\coloneqq \|\g(s)-\g(s')\|.
$$
Then, one checks easily
\beq{gen005}
\dpr_s d(s,s')=-\cos \th,\qquad \dpr_{s'} d(s,s')=\cos \th'.
\eeq
In particular, if we lift everything to the universal cover and introduce the new
coordinates $(x,y)=(s,\cos \th)\in\real\times [-1,1]$, then the billiard map is a twist map with $d$ as generating function, and it preserves the area form $dx\wedge dy$.\footnote{See \cite{tabachnikov2005geometry, siburg2004principle} for details.} Let $\G\coloneqq \{x_n,y_n\}_{n\ge 0}\subset \real\times [-1,1]$ be a billiard trajectory. 
\dfn{rotnum}
If the limit $\lim_{n\to\pm\infty} x_n/n$ exists then it is called the rotation number of the trajectory $\G=\{x_n,y_n\}_{n\ge 0}$ and will be denoted by $\rot(\G)$.
\edfn

\rem{perirot}
The trajectory $\G\coloneqq \{x_n,y_n\}_{n\ge 0}$ is called $(p,q)$--periodic if $y_q=y_0$ and $x_q=x_0+p$.  For such orbits, the rotation number always exists and is equal to $p/q$. Furthermore, there always exists at least two periodic orbits of rotation number $p/q$, for any rational $p/q$ by a Theorem by Birkhoff \cite{birkhoff1913proof}.
\erem

\dfn{mathbetfun}
Given $p/q\in\rational$, denote by $\mathscr{P}er(p,q)$ the collection of all the $(p,q)$--periodic trajectories. Thus, the Mather's $\beta$--function at $p/q$ is defined as
$$
\b(\frac{p}{q})\coloneqq -\max\{ \len(\G)\ : \ \G\in\mathscr{P}er(p,q)\}.
$$
Then, $\b$ is extended continuously to $\real$ by setting
$$
\b(\o)\coloneqq \lim_{n\to\infty} \b(\frac{p_n}{q_n}),
$$
for any sequence $\{{p_n}/{q_n}\}_n$ of rationals converging to $\o$.
\edfn

\dfn{def02}
For a given domain $\O$,  we define the "normalized" Mather's $\beta$--function as the function $\l^{-3}(\b+\ell\,\id)$, where $\b$ is the Mather's $\beta$--function, $\l$ is the Lazutkin perimeter, and $\ell$ the perimeter of the domain $\O$. 
\edfn
 The Taylor coefficients of the Mather's $\beta$--function at $0$ are related to the so-called Marvizi--Melrose invariants which was introduced by S. Marvizi and R. Melrose in \cite{marvizi1982spectral} where they proved the following. For, given $q\in\natural$, denote by
 \begin{align*}
   L_q&\coloneqq \sup\{ \len(\G)\ : \ \G\in\mathscr{P}er(p,q), \ p\in \integer\}, \\
   l_q&\coloneqq \inf\{ \len(\G)\ : \ \G\in\mathscr{P}er(p,q), \ p\in \integer\},
 \end{align*}
 the supremum and
the infimum of the perimeters of simple periodic billiard trajectories with $q$ vertices, respectively. The:
 
 \thm{MarvMel0}For any positive integer $k$ we have 
 $$
 \lim_{q\to\infty} q^k(L_q-l_q)=0.
 $$
 Moreover, $L_q$ has an asymptotic expansion as $q\to\infty$:
 $$
 L_q \sim \ell_0+\sum_{k=1}^\infty \frac{\ell_k}{q^{2k}},
 $$
 where $\ell_0$ is the length of the billiard table and $\ell_k$ are constants, depending on the curvature of the table.
 \ethm
 This collection $\{\ell_k\}_{k\ge0}$ constitutes the Marvizi--Melrose spectral invariants (also called Marvizi--Melrose invariants).
 
\ \\
\noi
In the case of non--integrables domains, we prove that the Taylor coefficients  of the "normalized" Mather's $\beta$--function are invariant under smooth conjugacies (see Theorem~\ref{teoMan}). As a {\bf conditional consequence}, under Conjecture~\ref{conj1}, the sequence given by the ratios of the respective Marvizi--Melrose invariants of such two $C^\infty$--conjugated domains is given by a geometric sequence (see Corollary~\ref{cor1MM}). In contrast, we prove that any two elliptic billiard maps are $C^0$--conjugated near their respective boundaries, and $C^\infty$--conjugated, in the open cylinder,  near their respective boundaries and off two lines. We also prove that the billiard maps corresponding to two ellipses are topologically conjugated only if the two ellipses are similar (\ie they are the same up to a rescaling and an isometry) or, equivalently, have the same eccentricity (see Theorem~\ref{teoMan2}).\footnote{It is an easy exercise to check that: two ellipses are similar \sss they have the same eccentricity.} Finally, assuming the Birkhoff conjecture, we prove  that the answer to Question~\ref{qes1} is affirmative whenever one of the domains is integrable (See Corollary~\ref{guillint}).
\section{Results}
\thm{teoMan} Let $\O_1,\O_2\subset \rn$ be two strictly convex domains with smooth boundaries. Let $\b_j$ be the Mather's $\b$--function of $\O_j$, $\ell_j\coloneqq \len(\dpr\O_j)$, $\kappa_j(l)$ be the curvature of $\dpr\O_j$ with arc--length parametrization $l$ and $\l_j\coloneqq \int_0^{\ell_j}\kappa_j(l)^{2/3}dl$.\footnote{Observe that $\kappa_j>0$ by strictly convexity and, therefore, $\l_j>0$. 
} Assume that no strip in $\O_2\times [0,\pi)$ containing the boundary $\O_2\times\{0\}$ is foliated by invariant curves 
and
\beq{conjf1f2}
f_{\O_1}\circ \tilde{h}=\tilde{h}\circ f_{\O_2},
\eeq
for some homeomorphism $\tilde{h}\in C^\infty(\torus\times\real)$. Then, $\l_1^{-3}(\b_1(\o)+\ell_1\o)$ and $\l_2^{-3}(\b_2(\o)+\ell_2\o)$ have the same Taylor's expansion at $\o=0$. In particular, $\l_1^{-3}\mathscr L_1(\o)$ and $\l_2^{-3}\mathscr L_2(\o)$ have the same Taylor's expansion at $\o=0$, where $\mathscr L_j(\o)$ is the Lazutkin's parameter of $\O_j$ corresponding to the convex caustic of rotation number $\o$ whenever it exists.
\qed
\ethm
Denote by by $\{\mathcal{I}^{j}_{2n+1}\}_n$ the Marvizi--Melrose invariants of $\O_j$. We conjecture the following relation between the Taylor's coefficients of the Mather's $\b$--function and the Marvizi--Melrose invariants.\footnote{The Marvizi--Melrose invariants are indeed known to be algebraically equivalent to the Taylor's coefficients of the Mather's $\b$--function but, to our best knowledge, no explicit algebraic formula relating them in known.}

\conj{conj1} For any $n\in \natural$,
\beq{betaMMRe}
\frac{d^{2n+1}\b}{d\o^{2n+1}}(0)=\mathcal{I}_3^{n+2}\cdot\sum_{\substack{\s_1,\cdots,\s_n\in\natural\cup\{0\}\\ \s_1+\cdots+\s_{n}=n-1\\ 3\s_1+\cdots+(2n+1)\s_{n}=5(n-1)}}r_n(\s_1,\cdots,\s_n)\cdot \mathcal{I}_3^{\s_1}\cdots \mathcal{I}_{2n+1}^{\s_{n}}\;,
\eeq
where each $r_n(\s_1,\cdots,\s_n)\in\rational$. \qed
\econj
\noi
Under Conjecture~\ref{conj1}, Theorem~\ref{teoMan} yields easily (\eg by induction) the following.
\cor{cor1MM}
The sequence $\{\mathcal{I}^{1}_{2n+1}/\mathcal{I}^{2}_{2n+1}\}_n$ is a geometric sequence:\footnote{Observe that $\l_j=\mathcal{I}^j_3$.}
\beq{RelMM}
\mathcal{I}^{1}_{2n+1}=(\l_2/\l_1)^{2n-3}\;\mathcal{I}^{2}_{2n+1},\qquad\forall\; n\ge 1.
\eeq
\ecor
\rem{conjBetf}
Using \cite{sorrentino2015computing}, one checks easily \equ{betaMMRe} and \equ{RelMM} up to $n=4$. 
\erem
In contrast to Theorem~\ref{teoMan}, the following holds for elliptic Billiards and, hence, according to Birkhoff Conjecture~\ref{conj1}, for convex integrable billiards.
\thm{teoMan2}
Let $\O_j$ be an ellipse with eccentricity $\mathbf{e}_j$ and $f_{j}$ be the associated billiard map ($j=1,2$). Set\footnote{We use $-p$ to denote the symmetric of the point $p$ \wrt the origin.} 
$$\O_j^*\coloneqq \O_j\setminus\{p_j^*,-p_j^*\}, \quad\mbox{where}\quad p_j^*\coloneqq (a_j\cos \arctan^{-1}(1-\mathbf{e}_j^2),-b_j\sin \arctan^{-1}(1-\mathbf{e}_j^2) ).$$ 
Then

\noi
$(i)$ The billiard maps $f_1$ and $f_2$ are $C^0$--conjugated near their respective boundaries {\it i.e.}, there exists $\th_j^*>0$ such that their restrictions $f_1|_{\O_1\times (0,\th_1^*)}$ and $f_2|_{\O_2\times (0,\th_2^*)}$ are $C^0$--conjugated. Moreover,  the restrictions $f_1|_{\O_1^*\times (0,\th_1^*)}$ and $f_2|_{\O_2^*\times (0,\th_2^*)}$ are $C^\infty$--conjugated.

\noi
$(ii)$ If the billiard maps $f_1$ and $f_2$ are (globally) $C^0$--conjugated then $\O_1$ and $\O_2$ are similar (or, equivalently, $\mathbf{e}_1=\mathbf{e}_2$).
\ethm

\noi
As a consequence, assuming the Birkhoff conjecture, we have: 
\cor{guillint}
 Let $f$ and $g$ be smooth Birkhoff billiard maps corresponding to two
strictly convex domains $\O_f$ and $\O_g$. Assume that $f$ and $g$ are topologically conjugated and that one of the two domains is integrable. Then, $\O_f$ and $\O_g$ are similar. 
\ecor
\section{Proofs}
\subsection{Proof of Theorem~\ref{teoMan}\label{parteo1}}
 Given $j=1,2$, pick a symplectic change of coordinates $\phi_j$  with associated Cantor set $\mathscr C_j$  corresponding to $f_{\O_j}$ given by  Theorem~\ref{KovPop}; set
$$
f_j\coloneqq \phi_j^{-1}\circ f_{\O_j}\circ \phi_j,\quad \mbox{and}\quad h=(h_1,h_2)\coloneqq \phi_1^{-1}\circ \tilde{h}\circ \phi_2.
$$
Let $\b_j\coloneqq \b_{\O_j}$, the minimal average action (or Mather's $\b$--function) of $f_{\O_j}$. 
\noi
Observe that \equ{conjf1f2} reads
\beq{conjf1f2bis}
f_{1}\circ h= h\circ f_{2},
\eeq
which, recalling\footnote{We refer the reader to \cite{sorrentino2015computing} for explicit formula of the coefficients of the Taylor's expansion of the $\a$ and $\b$ functions at the boundary.} $\a_1'(-\ell_1)=\a_1'(-\ell_1)=0$, implies for $\th=-\ell_2$,
$$
\a_1'(h_2(s,-\ell_2))=0=\a_1'(-\ell_1),\qquad \forall\; s\in \torus.
$$
Thus,
\beq{heell2}
h_2(s,-\ell_2)=-\ell_1,\qquad \forall\; s\in \torus.
\eeq
Given $j=1,2$ and $\o\in\a_j'(\mathscr C_j)$, let denote by $\G^j_\o\coloneqq \torus\times \{\b_j'(\o)\}$, the KAM torus of frequency $\o$ of $f_j$ and by $\mathscr D^j_\o$, the domain of the cylinder $\torus\times\real$ bounded by the closed curves $\{\th=\b_j(\o)\}$ and $\{\th=-\ell_j\}$.

\noi
Combining \equ{conjf1f2bis} and \equ{conjBill}, we get that, for any $\o\in\a_2'(\mathscr C_2)\setminus\{0\}$, $h(\G^2_\o)$ is an invariant curve of Diophantine frequency $\o$ for $f_1$. Thus, by Lemma~\ref{LagTor}, $h(\G^2_\o)$ is Lagrangian and, hence, is a KAM curve of frequency $\o$ of $f_1$. Therefore, by uniqueness, for any $\o\in\a_2'(\mathscr C_2)\setminus\{0\}$, 
\beq{UniKAm}
h(\G^2_\o)=\G^1_\o,
\eeq
and, by continuity, for any $(s,\th)\in\torus\times\mathscr C_2$,
\beq{Alf1Alf2}
\a_1'(h_2(s,\th))= \a_2'(\th).
\eeq
Thus, as $h$ is a homeomorphism of the cylinder $\torus\times\real$, we infer
\beq{boudry}
h(\mathscr D^2_\o)=\mathscr D^1_\o\qquad \mbox{and}\qquad h(\dpr\mathscr D^2_\o)=\dpr\mathscr D^1_\o.
\eeq
Set $\jac h\coloneqq \det \dpr_{(s,\th)} h$.\footnote{ Here and thereafter, we denote by $\dpr_{(s,\th)} h$ the Jacobian matrix of $h$.} We claim that
\beq{dthhOinfo}
\dpr_\th^k\jac h(s,-\ell_2)=0,\qquad\forall\; s\in\torus,\ \forall\; k\ge 1.
\eeq
Fix $s\in\torus$. Indeed, we are going to construct inductively, for a given $k\ge 1$, a sequence $\{\o_n^k\}_n\searrow 0$ such that 
\beq{dthhOinfoint}
\dpr_\th^k\jac h(s,\b_2'(\o_n^k))=0, \qquad\forall\; n\ge 1,
\eeq
which implies \equ{dthhOinfo} by continuity. For, pick a strictly decreasing sequence $\{\o_n^0\}_n\subset \a_2'(\mathscr C_2)\setminus\{0\}$ converging to $0$.\footnote{Such sequence exists as $\a_2'(\mathscr C_2)$ is a perfect set.}
As no strip near the boundary of $\O_2$ is foliated by invariant curves,   
then, by Mather connecting Theorem (see \cite{mather1990differentiability,kaloshin2003geometric}) any two successive invariant curves  of $f_2$ are (asymptotically) connected by some orbit, the invariant curves being naturally ordered via their respective rotation numbers. Let then $\mathscr M_2$ be the closure of the union of all invariant curves  and Mather's connecting orbits of $f_2$; in particular, $\mathscr M_2$ contains the family of KAM tori of $f_2$. As $f_1$ and $f_2$ are symplectic, we have $\det \dpr_{(s,\th)} f_1=\det \dpr_{(s,\th)} f_2\equiv 1$ and, therefore, differentiating \equ{conjf1f2bis}, we obtain, for any $(s,\th)\in \rn$, $\jac h(f_2(s,\th))= \jac h(s,\th)$ \ie $\jac h$ is constant along orbits of $f_2$. Consequently, $\jac h$ is constant on $\mathscr M_2$ by continuity. 
Then, for any $n\in \natural$, as $\jac h(s,\b_2'({\o}_n^0))=\jac h(s,\b_2'({\o}_{n+1}^0))$, by the Mean Value Theorem,  there exists $\o_{n+1}^0< {\o}_n^1< \o_n^0$ such that
$
\dpr_\th \jac h(s,\b_2'({\o}_n^1))=0.
$
Observe that, by construction, $\{\o_n^1\}_n\searrow0$, which proves \equ{dthhOinfoint} for $k=1$. Moreover, we recover the setting we started with, with now $\{\o_n^1\}_n$ playing the role of $\{\o_n^0\}_n$. Iterating the argument, we then obtain \equ{dthhOinfoint} and, in particular, \equ{dthhOinfo}. Hence, for $\o$ near $0$,\footnote{Recall that $\b_2'(\o)+\ell_2=O(\o^2)$; see \cite{sorrentino2015computing}.}
\beq{oDInf}
\jac h(s,\b_2'(\o))\eqby{dthhOinfo} \jac h(s,-\ell_2)+O((\b_2'(\o)+\ell_2)^k)= \jac h(s,-\ell_2)+O(\o^{2k}),\qquad \forall\; k\ge 1.
\eeq


\noi
Thus, for any $\o\in \a_1'(\mathscr C_1)$ and $k\ge 1$,
\begin{align}
\b_1'(\o)+\ell_1=\int_{\dpr \mathscr D_\o^1}\th ds&= \int_{\mathscr D_\o^1}d\th ds \qquad \mbox{(by Stokes--Whitney's Theorem)}\nonumber\\
  &= \int_{h(\mathscr D_\o^1)}|\jac h|d\th ds\eqby{boudry} \int_{\mathscr D_\o^2}|\jac h|d\th ds \quad \mbox{(by making $(s',\th')=h(s,\th)$)}\nonumber\\
  &=\int_{\dpr\mathscr D_\o^2}|\jac h|\th ds\qquad\qquad\qquad\qquad\quad \mbox{(by Stokes--Whitney's Theorem)}\nonumber\\
  &= \int_{\torus}(|\jac h(s,\b'_2(\o))|\b'_2(\o)+|\jac h(s,-\ell_2)|\ell_2)\; ds\nonumber\\
  &\eqby{oDInf} \int_{\torus}((\b_2'(\o)+\ell_2)|\jac h(s,-\ell_2)|+O(\o^{2k}))\; ds\nonumber\\
  &= (\b_2'(\o)+\ell_2)\int_{\torus}|\jac h(s,-\ell_2)| ds+O(\o^{2k})\label{impReb1b2}
\end{align}
We now compute $\int_{\torus}|\jac h(s,-\ell_2)| ds$. In fact, by Lemma~\ref{lem1}, we have $\jac h(s,\th)\equiv h_2'(\th)$. But\footnote{See \cite{sorrentino2015computing} .}
$$
(\th+\ell_j)^{-1/2}\a_j'(\th)= 2\sqrt{2}\l_j^{-3/2}+ O(\th+\ell_j).
$$
Thus, rewriting the second part of \equ{InThbis} 
 as
$$
\left(\frac{h_2(\th)+\ell_1}{\th+\ell_2}\right)^{1/2}(h_2(\th)+\ell_1)^{-1/2}\a_1'(h_2(\th))= (\th+\ell_2)^{-1/2}\a_2'(\th),
$$
and using \equ{heell2}, it follows
$$
\jac h(s,-\ell_2)\equiv h_2'(-\ell_2)= \left(\frac{\l_1}{\l_2}\right)^3.
$$
Hence, for any $k\ge 1$,
$$
\l_1^{-3}(\b_1'(\o)+\ell_1)=\l_2^{-3}(\b_2'(\o)+\ell_2)+O(\o^{2k}),
$$
and, therefore,
\beq{eqAux1}
\l_1^{-3}(\b_1(\o)+\ell_1\o)=\l_2^{-3}(\b_2(\o)+\ell_2\o)+O(\o^{2k+1}),
\eeq
\ie $\l_1^{-3}(\b_1(\o)+\ell_1\o)$ and $\l_2^{-3}(\b_2(\o)+\ell_2\o)$ have the same Taylor expansion. Now, for $\o$ small enough, if a convex caustic of rotation number $\o$ of $\O_j$ exists, then it Lazutkin's parameter  is given by  $\mathscr L_j(\o)=\o\b_j'(\o)-\b_j(\o)$ and, therefore, $\l_1^{-3}\mathscr L_1(\o)=\l_2^{-3}\mathscr L_2(\o)+O(\o^{2k+1})$.
\qed

\lem{lem1}
Under the notations and assumptions above, there exists\footnote{Actually, $g_1$ is as smooth as $h$.} $g_1\in C^\infty(\real,\real)$ such that, for any $(s,\th)\in\torus\times\mathscr C_2$, with $\th$ sufficiently close to $-\ell_2$,
\beq{InThbis}
h(s,\th)=(g_1(\th)+ s,h_2(\th)),\qquad \a_1'(h_2(\th))= \a_2'(\th).
\eeq
\elem

	\proof
Let $(s,\th)\in \torus\times \mathscr C_2\setminus\{-\ell_2\}$. Then, we have
\begin{align}
\big(h_1(s,\th)+\a_1'(h_2(s,\th)\big), h_2(s,\th))\eqby{conjBill} f_1\circ h(s,\th)\eqby{conjf1f2} h\circ f_2(s,\th)\eqby{conjBill}\nonumber\\
\qquad\eqby{conjBill}\big(h_1(s+\a_2'(\th),\th),h_2(s+\a_2'(\th),\th) \big).\label{deteH}
\end{align}
Hence, $h_2(s,\th)=h_2(s+\a_2'(\th),\th)$ which implies $h_2$ is independent of $s$ \ie $h_2(s,\th)=h_2(\th)$, as $\a_2'(\th)\not\in \rational$. We have also
\beq{h1o}
h_1(s,\th)+\a_1'(h_2(\th))\eqby{deteH} h_1(s+\a_2'(\th),\th),
\eeq
which differentiated \wrt $s$ yields $\dpr_s h_1(s,\th)= \dpr_s h_1(s+\a_2'(\th),\th)$. Thus, as $\a_2'(\th)\not\in \rational$, $\dpr_s h_1$ is also independent of $s$: $\dpr_s h_1(s,\th)=g_2(\th)$,  so that $h_1(s,\th)=g_1(\th)+g_2(\th)\cdot s$, for some $g_1,g_2\in C^\infty(\real,\real)$. Now, plugging the expression of $h_1$ found into \equ{h1o} and using \equ{Dalf:bet}, we get 
\beq{InTh}
h(s,\th)=(g_1(\th)+g_2(\th)\cdot s,h_2(\th)),\qquad \a_1'(h_2(\th))= g_2(\th)\a_2'(\th).
\eeq
Now, combining \equ{Alf1Alf2},  \equ{InTh} with the fact that\footnote{See \cite{sorrentino2015computing}.} $\a_2'(\th)\not=0$ for $\th+\ell_2>0$ small enough yields
$
g_2(\th)=1,
$
for $\th$ near $-\ell_2$. Thus, as $\mathscr C_2\setminus\{-\ell_2\}$ accumulates at $-\ell_2$, by continuity we obtain $g_2\equiv 1$, and \equ{InThbis} is proven. 

\qed
\rem{remteoMan2}
{\bf (i)}
Observe that, according to the proof, the conclusion in Lemma~\ref{lem1} still holds if $\tilde{h}$ is assumed merely to be $C^2(\torus\times\real,\real)$.\\ 
{\bf (ii)} Observe also that, if $\O_1$ and $\O_2$ are both ellipses and $\tilde{h}$ is $C^2$ conjugacy near their respective boundaries of the associated billiard maps $f_{\O_1}$ and $f_{\O_2}$, then \equ{UniKAm} and \equ{Alf1Alf2}, and, therefore,  \equ{InThbis} hold is some neighborhood $[-\ell_2,-\ell_2+\vae_2)$ of $-\ell_2$. Indeed, by the same argument, one gets that \equ{UniKAm}, \equ{Alf1Alf2} and  \equ{InThbis} hold on $[-\ell_2,-\ell_2+\vae_2)\setminus\rational$ and, therefore, by continuity they hold on $[-\ell_2,-\ell_2+\vae_2)$.
\erem
\subsection{Proof of   Theorem~\ref{teoMan2}}
\subsubsection{Proof of $(i)$ in  Theorem~\ref{teoMan2}\label{locconj}}
We shall exhibit a $C^\infty$--conjugacy. Given $j=1,2$ let $0<b_j\le a_j$ such that
$$
\O_j\coloneqq \left\{(x,y)\in\rn\ : \ \frac{x^2}{a_j^2}+\frac{y^2}{b_j^2}=1\right\},
$$
 so that $\mathbf{e}_j\coloneqq \sqrt{1-(b_j/a_j)^2}$. 
 Let  $\mathsf{c}_j\coloneqq \sqrt{a_j^2-b_j^2}$ its semi--focal distance, and $\mathcal{F}_j^\pm\coloneqq (\pm \mathsf{c}_j,0)$ its two foci. Consider the family of confocal elliptic caustics 
 $$
\O_{j,\l}\coloneqq \left\{(x,y)\in\rn\ : \ \frac{x^2}{a_j^2-\l^2}+\frac{y^2}{b_j^2-\l^2}=1\right\}\;,\qquad \l\in[0,b).
$$
We denote the (incomplete) elliptic integral of the first kind by
$$
F\colon [0,2\pi)\times[0,1)\ni (\f,k)\longmapsto \int_0^\f\frac{d\phi}{\sqrt{1-k^2\sin^2\phi}},
$$
and the complete elliptic integral of the first kind by $K(k)\coloneqq F(\pi/2,k)$.

\noi
 Consider the following auxiliary functions 
 \begin{align*}
 &k_{j}\colon [0,b_j)\ni\l\longmapsto \sqrt{\frac{a_j^2-b_j^2}{a_j^2-\l^2}}\in[0,1)\;,\qquad \o_j\colon [0,b_j)\ni\l\longmapsto\frac{F(\arcsin(\l/b_j),k_j(\l))}{2K(k_j(\l))}\in [0,\su2)\;,\\
 &h_{j,1}\colon [0,2\pi)\times [0,\pi)\ni (\phi,\th)\longmapsto (a_j\cos\phi,b_j\sin\phi,\th)\in \O_j\times [0,\pi),\\
  &h_{j,2}\colon [0,2\pi)\times [0,\pi)\ni (\phi,\th)\longmapsto (\l,t)\coloneqq \bigg({g_j(\phi,\th)}, F\big(\phi,k_j({g_j(\phi,\th)}\big)\bigg)\in [0,b_j)\times [0,4K(k_j(b_j)))\\
  &h_{j,3}\colon [0,1)\times [0,1)\ni (\l,t)\longmapsto \big(b_j\l,4K(k_j(b_j\l))t\big)\in [0,b_j)\times [0,4K(k_j(b_j)))\;,\\  
 \end{align*}
 where
\beq{gjdef}
 \begin{aligned}
 g_{j}(\phi,\th)\coloneqq& 
 \frac{|a_j^2\sin\phi+b_j^2\cos\phi|\sin\th}{\sqrt{a_j^2\sin^2\phi+b_j^2\cos^2\phi}}\;.
 \end{aligned}
 \eeq
It is obvious that each of the maps $k_j$, $\o_j$, $h_{j,1}$ and $h_{j,3}$ are $C^\infty$--diffeomorphism (onto their respective codomains). The map $h_{j,2}$ is a homeomorphism (onto its codomain) and its restriction $h_{j,2}^*\colon \mathscr T^*\times (0,\pi)\overset{onto}\longrightarrow (0,b_j)\times [0,4K(k_j(b_j)))$ is a  $C^\infty$--diffeomorphism, where $\mathscr T^*\coloneqq [0,2\pi)\setminus \{-\arctan^{-1}(b_j^2/a_j^2),\pi-\arctan^{-1}(b_j^2/a_j^2)\}$.\footnote{Indeed, the set of singularities of $h_{j,2}$ is precisely $g_j^{-1}(\{0\})=\{(\phi,\th)\in[0,2\pi)\times [0,\pi):\ \tan\phi=-b_j^2/a_j^2\mbox{ or } \th=0\}$} In particular, the following map is a $C^\infty$--diffeomorphism as well:
$$
h_{j,4}\colon [0,1)\times[0,1)\ni (\l,t)\longmapsto (\o_1^{-1}\circ \o_2(\l),t)\in [0,1)\times[0,1).
$$
Given $s=(x,y)\in\O_j$, let\footnote{Here and thereafter, $[P,Q]$ denotes the closed segment joining $P$ and $Q$ \ie $[P,Q]\coloneqq \{\a P+(1-\a)Q\;:\; 0\le \a\le 1\}$.}
\begin{align*}
\th_j(s)&\coloneqq \sup\left\{\vae\in (0,\pi)\ :\ \bigcup_{n=0}^\infty [f_j^n(s,\th),\;f_j^{n+1}(s,\th)]\bigcap [\mathcal{F}_j^-,\;\mathcal{F}_j^+]=\emptyset\;,\ \forall\; \th\in [0,\vae]\right\}\;.
\end{align*}
Then, one checks easily that
$$
\th_j(x,0)=\frac{\pi}{2}\qquad \mbox{and} \qquad \th_j(x,y)=\arctan\frac{b_j^2}{\mathsf{c}_j|y|},\quad \mbox{for}\ y\not=0.
$$
In particular, 
$$
\th_j^*\coloneqq \inf_{s\in\O_j}\th_j(s)= \arctan\frac{b_j^2}{\mathsf{c}_j}>0.
$$
Let  $\th_3^* \in [0,\pi)$ such that $h^{-1}(\O_1\times[0,\th_1^*))=\O_2\times[0,\th_3^*)$ and define 
$$
\th^*\coloneqq \min\{\th_2^*,\th_3^*\}.
$$
Furthermore, consider 
$$
h_j\coloneqq h_{j,1}\circ h_{j,2}^{-1}\circ h_{j,3}\circ h_{j,4}:[0,1)\times [0,1)\to \O_j\times [0,2\pi)
$$
and 
$$
h\coloneqq h_1\circ h_2^{-1} \colon\O_2\times [0,\pi)\to \O_1\times [0,\pi).
$$
Then, $h$ is a homeomorphism and its restriction
$$
h^*\colon\O_2^*\times (0,\pi)\to \O_1^*\times (0,\pi).
$$
is a $C^\infty$--diffeomorphisms. 
%
%
%
\noi
Now, we claim that
\beq{MagCoj}
f_1\circ h= h\circ f_2 \qquad\mbox{on}\qquad \O_2\times[0,\th^*),
\eeq
which would prove $(i)$ in Theorem~\ref{teoMan2}. The prove of \equ{MagCoj} follows easily a result proven in \cite{chang1988elliptical} (see also {\cite[Proposition~16]{kaloshin2018local}}), which can be reformulated as follows. The main difference with the aforementioned result is that we provide the explicit value of $(\l,t)$, the ``action--angle'' coordinates (see Appendix~\ref{applamb} for the computation of $\l$).  
\lem{lemellp1}
Let $(s,\th)\in \O_j\times [0,\pi)$. Set
$
(\phi,\th)\coloneqq h_{j,1}^{-1}(s,\th),\ (\l,t)\coloneqq h_{j,2}(\phi,\th)$ and  $\d_{j}(\l)\coloneqq 2F(\arcsin(\l/b_j),k_j(\l))=4K(k_j(\l))\o_j(\l).
$ 
Then, the segment joining $(s,\th)$ and $f_j(s,\th)$ is tangent to $\O_{j,\l}$ and
\beq{caustelip}
f_j\circ h_{j,1}\circ h_{j,2}^{-1}(\l,t)=f_j(s,\th)= h_{j,1}\circ h_{j,2}^{-1}(\l,t+\d_j(\l)).
\eeq
\elem

\subsubsection{Proof of $(ii)$ in  Theorem~\ref{teoMan2}\label{locconjglo}}
We adopt the notations in \S \ref{locconj}. Given $j=1,2$ and two coprime integers $0<m<n/2$, consider the function
$$
g_j\colon (-\mathsf{c}_j^2,0)\ni \x\longmapsto F\left(\arcsin\sqrt{\frac{b_j^2}{b_j^2-\x}},\sqrt{1+\frac{\x}{c_j^2}} \right)-\frac{2m}{n}F\left(\frac{\pi}{2},\sqrt{1+\frac{\x}{c_j^2}} \right).
$$
Then\footnote{See \cite[pg. 4569, Eq. (22)]{sieber1997semiclassical}.}, $f_j$ has a periodic orbit of rotation $m/n$ 
whose caustic is a hyperbola \sss the equation
\beq{hypPer}
g_j(\x)=0,\qquad \x\in (-\mathsf{c}_j^2,0)
\eeq
has a solution. But, we have
$$
g_j(-\mathsf{c}_j^2)=\arcsin\frac{b}{a}-\frac{m}{n}\pi,\qquad \mbox{and}\qquad \lim_{\x\to 0^-} g_j(\x)= +\infty.
$$
Hence, if 
\beq{condPe1}
\frac{m}{n}\ge \frac{1}{\pi}\arcsin\frac{b_j}{a_j},
\eeq
then \equ{hypPer} has a solution. 
{ In contrast, if 
\beq{condPe2}
\frac{m}{n}< \frac{1}{\pi}\arcsin\frac{b_j}{a_j},
\eeq
then \equ{hypPer} has no solution. Indeed, assuming \equ{condPe2}, we have, for any $\x\in (-\mathsf{c}_j^2,0)$,
$$
g_j(\x)> 
u_j(\x)>0,
$$
where
$$
u_j(\x)\coloneqq F\left(\arcsin\sqrt{\frac{b_j^2}{b_j^2-\x}},\sqrt{1+\frac{\x}{c_j^2}} \right)-\frac{2}{\pi}\arcsin(\frac{b_j}{a_j})F\left(\frac{\pi}{2},\sqrt{1+\frac{\x}{c_j^2}} \right).
$$
and, using \textit{Wolfram Matematica}\cite{Mathematica}, one gets
$$
\min\left\{u_j(\x):\ a_j> b_j>0,\ -a_j^2+b_j^2<\x<0 \right\}=5.65246\cdot 10^{-9}.
$$
}
Now, assume $f_1$ and $f_2$ are (globally) $C^0$--conjugated. By contradiction, assume $\mathsf{e}_1>\mathsf{e}_2$ \ie ${b_1}/{a_1}<{b_2}/{a_2}$. Then, for any two coprime integers $0<m<n/2$ such that,
$$
 \frac{1}{\pi}\arcsin\frac{b_1}{a_1}\le  \frac{m}{n}< \frac{1}{\pi}\arcsin\frac{b_2}{a_2}\,,
$$
$f_1$ admits a periodic orbit of rotation $m/n$ 
whose caustic is a hyperbola while $f_2$ does not admit such a periodic orbit, contradiction. Thus, $\mathsf{e}_1=\mathsf{e}_2.$
\qed

\section*{Appendix}
\addcontentsline{toc}{section}{Appendices}
\setcounter{section}{0}
\renewcommand{\thesection}{\Alph{section}}

\appA{Auxiliary facts about billiard dynamics}

\thmtwo{TeoAlf}{\cite{mather1994action, mather1990differentiability,siburg2004principle}}
Let $f$ be a monotone twist map. 
Then:
\begin{itemize}
\item[$(i)$] $\b$ is strictly convex. In particular, it is continuous and admits a Legendre--Fenchel transform 
$$
\a(c)=\sup_{\o\in\real} \o\cdot c-\b(\o),
$$
the so--called Mather's $\a$ function. Moreover, $\a$ is convex and we have 
\beq{Dalf:bet}
\a'(\b'(\o))=\o,
\eeq
at each point of differentiability $\o$ of $\b$.
\item[$(ii)$] $\b$ is differentiable at any irrational.
\item[$(iii)$] $\b$, and hence $\a$, is an invariant under any symplectic change of coordinates.
\end{itemize}
\ethm


\thmtwo{KovPop}{\cite{kovachev1990invariant, siburg2004principle}}
Let $\O\subset \rn$ be a smooth strictly convex closed curve of length $\ell$ and $f$ its corresponding billiard map. Then, there exists $\vae_0>0$ small, a cantor set $\mathscr C\subset [-\ell,-\ell+\vae_0)$, with $-\ell\in\mathscr C$, a smooth (exact) symplectomorphism $\phi\colon \torus\times [-\ell,-\ell+\vae_0)\ni (x_0,\th_0)\longmapsto (x_1,\th_1)\in \phi(\torus\times [-\ell,-\ell+\vae_0))\subset \torus\times \real$ such that 
\beq{conjBill}
\phi^{-1}\circ f\circ \phi (x,\th)=(x+\a'(\th),\th),\qquad \forall \; (x,\th)\in \torus\times \mathscr C.
\eeq
In particular, the KAM curve of rotation number $\o\in \a'(\mathscr C)$ of $\phi^{-1}\circ f\circ \phi$ is given by the graph of the constant function $\torus\ni x\longmapsto \b'(\o)$.
\ethm

\appB{Isotropicity Lemma\label{appC}}
\lemtwo{LagTor}{\cite{herman1989inegalites,broer2009quasi}}
Let $f\colon \mathcal{M}\coloneqq \real^m\times\torus^m \to f(\mathcal{M})\subset \mathcal{M}$ be a symplectic\footnote{\wrt the canonical symplectic form.} diffeomorphism of class $C^1$
and $\phi: \torus^m\to \mathcal{M}$, a $C^1$--embedding. Assume that 
$$
f\circ \phi\circ g(x)= \phi\circ g(x+t\o), \qquad \forall\; x\in \torus^m,
$$
for some $\o\in\real^m$ and, that $\o$ is rationally independent\footnote{\ie for any $k\in \zn\setminus\{0\}$, $\o\cdot k\not=0$.} if $m\ge 2$ and $\o\in\real\setminus \rational$ if $m=1$. Then, the torus $\phi(\torus^m)$ is isotropic \ie $\mathit{i}^*\varpi\equiv 0$, where $\mathit{i}\colon \phi(\torus^m)\hookrightarrow\mathcal{M}$ is the inclusion map.
\elem

\appC{Computation of the ``action--angle'' coordinates for elliptic billiards\label{applamb}}
It is well--note that the caustic parameter $\l$, together with its conjugate $t\coloneqq F\big(\phi,k_j(\sqrt{\l})\big)$ is an action--angle coordinates for the billiard in the ellipse $\O_j$.\footnote{See \eg \cite{chang1988elliptical} or \cite[Proposition~16]{kaloshin2018local}.} We are only left to show that, $\l=g_j(\phi,\th)$ with $g_j$ defined as in \equ{gjdef} is indeed the parameter of the caustic corresponding to the billiard trajectory $\{f_j^{n}(h_{j,1}(\phi,\th)):\; n\in\integer\}$ starting at $h_{j,1}(\phi,\th)$. Let\footnote{For the sake of simplicity, we will drop the subscripts $j$.} $q_\l(t)=(q_\l^x(t),q_\l^y(t))=(a \ecn(t;k_\l), b \esn(t;k_\l))=(a\cos\phi,b\sin\phi)$, where $\ecn$ and $\esn$ are the Jacobi elliptic functions:
$$
\ecn(t;k_\l)\coloneqq \cos(\am (t;k_\l))\qquad\mbox{and}\qquad \esn(t;k_\l)\coloneqq \sin (\am (t;k_\l)),
$$
where $\am(t;k)$, called the amplitude of $t$, is given by: $t= F\big(\phi,k)\big)$ \sss $\phi=\am(t;k)$. Then,\footnote{see \cite{chang1988elliptical} or {\cite[Proposition~16]{kaloshin2018local}}.} the trajectory starting at $q_\l(t)$ and tangent to the caustic 
$$
\O_{\l}\coloneqq \left\{(x,y)\in\rn\ : \ \frac{x^2}{a^2-\l^2}+\frac{y^2}{b^2-\l^2}=1\right\}\;,\qquad \l\in[0,b),
$$
hits the boundary at $q_\l(t+\d_\l)$, where $\d_l\eqqcolon\d(\l)$ is given in Lemma~\ref{caustelip}. Then, the point of tangency $p_\l(t)$ of the segement $[q_\l(t),q_\l(t+\d_\l)]$ to the caustic $\O_\l$ is given by $p_\l(t)=(\sqrt{a^2-\l^2}\cos\f, \sqrt{b^2-\l^2}\sin\f)$, for some $\f\in [0,2\pi)$. Thus, the slope of the line $(q_\l(t),q_\l(t+\d_\l))$ is given by:
\begin{itemize}
\item[$\bullet$] $-\sqrt{\frac{b^2-\l^2}{a^2-\l^2}}\cot\f$, using its tangency property to $\O_\l$ at $p_\l(t)$;
\item[$\bullet$] $\frac{b\sin\phi-\sqrt{b^2-\l^2}\sin\f}{a\cos\phi-\sqrt{a^2-\l^2}\cos\f}$, using its points $q_\l(t)$ and $p_\l(t)$;
\item[$\bullet$] $\tan\left(\th-\arctan\frac{b\cos\phi}{a\sin\phi}\right)=\frac{a\tan\phi\tan\th-b}{a\tan\phi+b\tan\th}$, using the fact that the slope of the tangent to $\O$ at $q_\l(t)$ is $-\frac{b\cos\phi}{a\sin\phi}$ and that the angle between this tangent and the line $(q_\l(t),q_\l(t+\d_\l))$ is $\th$.
\end{itemize}
Hence,
\begin{align}
-\sqrt{\frac{b^2-\l^2}{a^2-\l^2}}\cot\f&=\frac{b\sin\phi-\sqrt{b^2-\l^2}\sin\f}{a\cos\phi-\sqrt{a^2-\l^2}\cos\f}\,,\label{slop1}\\
-\sqrt{\frac{b^2-\l^2}{a^2-\l^2}}\cot\f&=\frac{a\tan\phi\tan\th-b}{a\tan\phi+b\tan\th}\,.\label{slop2}
\end{align}
Now, \equ{slop1} is equivalent to:
\beq{slop3}
a\sqrt{b^2-\l^2}\cos\phi\cos\f+b\sqrt{a^2-\l^2}\sin\phi\sin\f=\sqrt{(a^2-\l^2)(b^2-\l^2)}\,.
\eeq
The relation  \equ{slop2} implies
\begin{align}
&\frac{\cos^2\f}{(a^2-\l^2)(a\tan\phi\tan\th-b)^2}=\frac{\sin^2\f}{(b^2-\l^2)(a\tan\phi+b\tan\th)^2}=\nonumber\\
&\qquad\qquad\qquad=\frac{\cos^2\f+\sin^2\f}{(a^2-\l^2)(a\tan\phi\tan\th-b)^2+(b^2-\l^2)(a\tan\phi+b\tan\th)^2}\nonumber\\
&\qquad\qquad\qquad=\frac{1}{(a^2-\l^2)(a\tan\phi\tan\th-b)^2+(b^2-\l^2)(a\tan\phi+b\tan\th)^2}\label{sloop1}
\end{align} 
so that
\begin{align}
&\frac{\cos\f}{\sqrt{a^2-\l^2}(a\tan\phi\tan\th-b)}\eqby{slop2}\frac{-\sin\f}{\sqrt{b^2-\l^2}(a\tan\phi+b\tan\th)}=\nonumber\\
&\qquad\quad\eqby{sloop1}\frac{\pm1}{\sqrt{(a^2-\l^2)(a\tan\phi\tan\th-b)^2+(b^2-\l^2)(a\tan\phi+b\tan\th)^2}}\,.\label{sloop11}
\end{align}
Thus, plugging into \equ{slop3} the expressions of $\cos\f$ and $\sin\f$ obtained from \equ{sloop11} and then squaring both sides of the relation obtained yields $\l=g(\phi,\th)$, which completes the computation.
\qed
\bibliographystyle{apa}
\bibliography{BibtexDatabase}
\end{document}